\documentclass 
{amsart}
\usepackage{amsmath}
\usepackage{amsthm}
\usepackage{amscd}
\usepackage{amsxtra}
\usepackage{amssymb}
\usepackage[all]{xy}
\usepackage{pb-diagram,pb-xy}  
\usepackage{boxedminipage}
\usepackage{ulem}
\usepackage{stmaryrd}
\usepackage[dvipdfm]{graphicx,color}

\newtheorem{theorem}{Theorem}[section]
\newtheorem{lemma}[theorem]{Lemma}
\newtheorem{proposition}[theorem]{Proposition}
\newtheorem{corollary}[theorem]{Corollary}
\newtheorem*{KOT}{Theorem of Kurokawa-Ochiai \cite{KO}}
\newtheorem*{MT}{Main Theorem}

\newtheorem*{H-S}{Theorem (Hovey-Sadofsky)}

\theoremstyle{definition}

\newtheorem{example}[theorem]{Example}

\newtheorem{def-pr}[theorem]{Definition-Proposition}
\newtheorem{def-th}[theorem]{Definition-Theorem}
\newtheorem{def-cor}[theorem]{Definition-Corollary}

\newcommand{\ord}{ \operatorname{ord} }

\newcommand{\Z}{{ \mathbb Z}}

\newcommand{\F}{{ \mathbb F}}

\newcommand{\N}{{ \mathbb N}}

\newcommand{\C}{{ \mathbb C}}

\newcommand{\dotbox}{\hbox to 1em{\hss.\hss}}

\begin{document}

\title{On the random variable $\N^r \ni (k_1, k_2, \ldots, k_r) \mapsto \gcd(n,k_1k_2\cdots k_r) \in \N$    }
\author{Norihiko Minami}
\address{Omohi College, Nagoya Institute of Technology, Gokiso, Showa-ku, Nagoya 466-8555}
\email{nori@nitech.ac.jp}

\maketitle

\section{Introduction and Main Theorem}

In his talk at JAMI 2009 on March 24, Kurokawa presented a rather mysterious looking identidity
of elementary number theory:
\begin{equation} \label{FI}
\frac{1}{n}\sum_{k=1}^n \gcd(n,k) = \prod_{p|n} \left(  1 + \left( 1 - \frac{1}{p} \right) \ord_p(n) \right)
\end{equation}

\begin{example}  $n=12 = 2^2\cdot 3$:
\begin{multline*}
\frac{1}{12} ( 1 + 2 + 3 + 4 + 1 + 6 + 1 + 4 +  3 + 2 + 1 + 12)  = \frac{40}{12} =  \frac{10}{3} \\
=  2\cdot  \frac{5}{3} = \left(  1 + \left( 1 - \frac{1}{2} \right)\cdot 2 \right)
\cdot \left(  1 + \left( 1 - \frac{1}{3} \right)\cdot 1 \right)
\end{multline*}
\end{example}

In fact, this is a special case (the case $r=1$) of identities obtained by \cite{KO}:


\begin{KOT}
For $n, r \in \N$,
\begin{equation} \label{KOI}
\frac{1}{n^r}\sum_{k_1, \ldots, k_r =1}^n \gcd(n,k_1\cdots k_r)
= 
\prod_{ p \mid n}  
\left[
\sum_{l = 0}^{r} {}_{ \ord_p(n) }H_l  \left( 1 - \frac{1}{p} \right)^l 
\right]
\end{equation}
Here, ${}_mH_r$ is the repeated combination:
\begin{align*}
{}_mH_r &:= \# \{ \text{degree $r$ homogeneous monomials in $m$ variables} \}  \\
&= {}_{r + (m-1)}C_{m-1} = {}_{m+r-1}C_r = (-1)^r \binom{-m}{r}  \\
\text{\it c.f.}\quad & ( 1 + x)^m = \sum_{r=0}^m {}_mC_rx^r,\quad (1-x)^{-m} \underset{|x|<1}=  
\sum_{r=0}^{\infty} \binom{-m}{r}(-x)^r = \sum_{r=0}^{\infty} {}_mH_r x^r
\end{align*}
\end{KOT}

\begin{example}  $n=6= 2\cdot 3,\ r=2$:
\begin{multline*}
\frac{1}{6^2} (  1 + 2 + 3 + 2 + 1 + 6 + 
2 + 2 + 6 + 2 +  2 + 6 +
3 + 6 + 3+ 6 + 3 + 6    \\
+2 + 2 + 6 + 2 + 2 + 6 +
1 + 2 + 3 + 2+ 1 + 6 +
6 + 6 + 6 + 6 + 6 + 6
)
\\
= \frac{133}{36} 
= \frac{7}{4}\cdot \frac{19}{9} = \frac{4+2+1}{4}\cdot \frac{9 + 6 + 4}{9}  \\
= \left(  1 + 1\cdot \left( 1 - \frac{1}{2} \right)  +  1\cdot \left( 1 - \frac{1}{2} \right)^2 \right)
\cdot \left(  1 + 1\cdot \left( 1 - \frac{1}{3} \right)  +  1\cdot \left( 1 - \frac{1}{3} \right)^2 \right)
\end{multline*}
\end{example}


\cite{KO} obtained \eqref{KOI} by studying some multivariable zeta function of Igusa type,
and Kurokawa said he is not aware of any elementary proof even for \eqref{FI}.

Now the purpose of this paper is to give a purely elementary proof of generalizations of
Kurokawa-Ochiai identities \eqref{KOI} from the view point of elementary probability theory.
Fixing $n\in\N$, we would like to understand the random variable:

\begin{equation} \label{old}
\begin{split}
\tilde{X}: \tilde{\Omega} := \N^r &\to \N   \\
       (k_1, k_2, \ldots, k_r) &\mapsto \gcd(n,k_1k_2\cdots k_r),
\end{split}
\end{equation}

Although $\tilde{\Omega} = \N^r$ is an infinite set, we would like to regard it being equipped with
the \lq\lq homogeneous measure\rq\rq.  For this purpose, we observe:
\begin{equation*}
k_i \equiv k'_i\ \mod n\   (i=1, 2, \ldots, r)\quad \implies\quad \gcd(n,k_1k_2\cdots k_r)  = \gcd(n,k'_1k'_2\cdots k'_r)
\end{equation*}
Thus, instead of \eqref{old}, we may equally consider the following random variable:
\begin{equation} \label{new}
\begin{split}
X: \Omega := \left\{ 1, 2, \ldots, n \right\}^r &\to \N   \\
       (k_1, k_2, \ldots, k_r) &\mapsto \gcd(n,k_1k_2\cdots k_r),
\end{split}
\end{equation}
where $\Omega= \left\{ 1, 2, \ldots, n \right\}^r$ is equipped with the homogeneous measure.
Then the Kurokawa-Ochiai identity \eqref{KOI} is nothing but a convenient formula to evaluate the average $E[X]$ of
the random variable \eqref{new}.   Now, from the view point of elementary probability theory, it is very natural to
seek for similar convenient formulae for the variance $V[X] = E[X^2] - E[X]^2$ and even \lq\lq higher\rq\rq invariants.  

Our Main Theorem 
offers such formulae for the continuous version $E[X^w]\ (w\in \C)$. Their
special cases $w\in\N$ are nothing but the moments of the randowm variable $X$, and the simplest case $w=1$
is nothing but the Kurokawa-Ochiai identity \eqref{KOI}.

\begin{MT} \label{MT}
For $n, r \in \N,\ w\in \C$,
{\small
\begin{align*}
&\quad \frac{1}{n^r}\sum_{k_1, \ldots, k_r =1}^n \gcd(n,k_1\cdots k_r)^w  \\
&=
\begin{cases}
\prod_{ p \mid n}
\left[
\left( \frac{ 1 - p^{-1} }{ 1 - p^{w-1} } \right)^r
+ p^{\ord_p(n)(w-1)} \sum_{l = 0}^{r-1} {}_{ \ord_p(n) }H_l  
\left\{
\left( 1 - p^{-1} \right)^{l} 
- 
\left( \frac{ 1 - p^{-1} }{ 1 - p^{w-1} } \right)^r
( 1 - p^{w-1} )^{l}   \right\} 
\right]
\
&\text{if}\ w\neq 1  \\
\prod_{ p \mid n}  
\left[
\sum_{l = 0}^{r} {}_{ \ord_p(n) }H_l  \left( 1 - \frac{1}{p} \right)^l 
\right]
\
&\text{if}\ w= 1
\end{cases}
\\
&=
\begin{cases}
\prod_{ p \mid n}
\left[
\left( \frac{ 1 - p^{-1} }{ 1 - p^{w-1} } \right)^r
+ p^{\ord_p(n)(w-1)} ( 1 - p^{-1} )^r \sum_{l = 0}^{r-1} {}_{ \ord_p(n) }H_l  
\left\{
\left( 1 - p^{-1} \right)^{l-r} 
- ( 1 - p^{w-1} )^{l-r}   \right\} 
\right]
\
&\text{if}\ w\neq 1  \\
\prod_{ p \mid n}  
\left[
\sum_{l = 0}^{r} {}_{ \ord_p(n) }H_l  \left( 1 - \frac{1}{p} \right)^l 
\right]
\
&\text{if}\ w= 1
\end{cases}
\end{align*}
}
\end{MT}

\begin{corollary}  \label{asymptotic} For $n\in \N$ and $w\in\C$ satisfying
$w=1$ 
or $ | p^w - 1 | > 1$,
\begin{equation} \label{limit}
\lim_{r\to \infty}  \frac{1}{n^r}\sum_{k_1, \ldots, k_r =1}^n \gcd(n,k_1\cdots k_r)^w  = n^w
\end{equation}
\end{corollary}

Heuristically speaking, \eqref{limit} implies the random variable \eqref{old} \lq\lq converges\rq\rq\ to the
constat random variable concentrated at $n$.

When $n$ is fixed, the number of the terms needed to be evaluaed in \newline
$\frac{1}{n^r}\sum_{k_1, \ldots, k_r =1}^n \gcd(n,k_1\cdots k_r)^w$ grow exponentially with respect to $r$, 
but Main Theorem allows us to evaluate it with lear growth with respect to $r$.

In our purely elementary proof for Main Theorem, the original identity \eqref{FI} turns out to be a matter of triviality.
In the course of proving the general case, we shall show some identitites
involving the repeated combinations (Lemma~\ref{HI}), which may of of independent interest.

We note another kind of generalizations of \eqref{FI} is given in \cite{DKK}, which
is once again proved using some zeta function of Igusa type.  Even for this and some generalization, we can
offer an elementary proof \cite{M1}.  
It is possible to define and study some multivarible {\it deformed} zeta function of $\F_1$-scheme of Igusa-type \cite{M2},
which generalizes both \cite{DKK} and \cite{KO}.

The basic idea of the results in this paper were obtained during the author's stay at JAMI2009, Johns Hopkins University, in March 2009,
and presented at NCGOA2009, Vanderbilt University, in May 2009.  The author would like to express his gratitude to 
Katia Concani, Alain Connes, Jack Morava, Takashi Ono, Steve Wilson, and Guoliang Yu for their hospitalities.
The author also would like to express his gratitude to Nobushige Kurokawa for his work and encouragement.

\section{Proof of Main Theorem}

Whereas [KO] used some multivariable Igusa type zeta funciions for 
groups,
we exploit some 
ring structure:
For a finite 
ring with $n$ elements $R = \{ k_i \in R \mid 1\leq i\leq n\}$ and $r\in \N$, set
\begin{equation*} 
Z_R^r
(w) := \frac{1}{ | R |^r } \sum_{  ( k_1, \ldots, k_r) \in R^r  } \big| R/ (k_1\cdots k_r)R 
\big|^w
\quad 
(w\in \C)
\end{equation*}
We wish to understand this, because
\begin{equation*}
Z_{ \Z/n\Z }^r
(w) =
\frac{1}{n^r}\sum_{k_1, \ldots, k_r =1}^n \gcd(n,k_1\cdots k_r)^{w}
\quad 
(w\in \C)
\end{equation*}

Of course, we have an elementary probability theoretical interpretation:
For the random variable
\begin{align*}
X: \Omega := R^r &\to \N   \\
       (k_1,k_2,\ldots,k_n) &\mapsto  \big| R/ (k_1{\color{green}{\cdots}} k_r)R \big|,
\end{align*}
where $\Omega = R^r$ is equipped with the homogeneous measure,
\begin{equation*}
E[X^w] =  \frac{1}{ | R |^r } \sum_{  ( k_1, \ldots, k_r) \in R^r  } \big| R/ (k_1{\color{green}{\cdots}} k_r)R 
\big|^{w}
= Z_R^r
(w)
\end{equation*}
\underline{\it Elementary properties of $Z_R^r(w)$:} 

\begin{itemize}
\item If there is a finite ring decomposition $R = \prod_{i} R_i$, then
\begin{equation*} \label{decomposition}
Z_R^r(w)
= \prod_i Z_{R_i}^r
(w)
\end{equation*}
\item Set
$
N_R^r(f) := \Big| \left\{ ( k_1, \ldots, k_r) \in R^r \mid | R/ (k_1\cdots k_r)R  | = f \right\} \Big|
$,
then
\begin{equation*}
Z_R^r(w)
:= \frac{1}{ | R |^r } \sum_{  ( k_1, \ldots, k_r) \in R^r  } 
\big| R/ (k_1\cdots 
k_r)R \big|^{w}
=  \frac{1}{ | R |^r } \sum_{f=1}^{ |R| } N_R^r(f)f^{w}
\end{equation*}
\end{itemize}

So, if $n = \prod_{ p | n } p^{\ord_p(n) }$, then by the Chinese Remainder Theorem,
\begin{equation} \label{reduction}
\begin{split}
\frac{1}{n^r}\sum_{k_1, \ldots, k_r =1}^n \gcd(n,k_1\cdots k_r)^{w}
&= Z_{ \Z/n\Z }^r(w) 
= \prod_{ p | n } Z_{ \Z/p^{\ord_p(n)}\Z }^r 
(w)  \\
&= \prod_{ p | n }  \frac{1}{ p^{r\ord_p(n)} } \sum_{f=1}^{ p^{\ord_p(n)} } N_{  \Z/ p^e\Z }^r(f) f^{w} 
\end{split}
\end{equation}
Thus, we need to compute 
\begin{align*}
N_{  \Z/ p^e\Z }^r(f) &= 
 \Big| \left\{ ( k_1, \ldots, k_r) \in (  \Z/ p^e\Z )^r \Big| | ( \Z/ p^e\Z )/ (k_1\cdots k_r)( \Z/ p^e\Z )  | = f \right\} \Big|      \\
&=  \Big| \left\{ ( k_1, \ldots, k_r) \in  \{ 1, 2, \ldots,  p^e \}^r \Big|  \gcd( p^e,  k_1\cdots k_r  )   = f \right\} \Big|   
\end{align*}
for $1\leq f\leq p^{\ord_p(n)}$.

To get some feeling, we first play with the simplest case $r=1$,
Then, it is easy to evaluate:
\begin{equation} \label{r=1_density}
N_{ \Z/ p^e\Z }^1( f )  =
\begin{cases}
p^{e-d}-p^{e-d-1}\quad &\text{if}\ f = p^d\ (0\leq d < e)   \\
1 \quad &\text{if}\ f = p^e \\
0 \qquad &\text{if otherwise}
\end{cases}
\end{equation}
which follows from
\begin{align*}
\# \Big| \{ 1, 2, \ldots,  p^e \} \Big| &= p^e = 
\left( \sum_{d=0}^{e-1}\underset{\ord_p=d}{\underbrace{(p^{e-d}-p^{e-d-1})}} \right) + 
\underset{\ord_p=e}{\underbrace{p^{e-e}}}
\end{align*}
%
%
%
\underline{\it Proof for $r=1$}:
\begin{align*}
&\quad
\frac{1}{n}\sum_{k =1}^n \gcd(n,k)^{w} 
\overset{\eqref{reduction}}{=} \prod_{ p | n }  \frac{1}{ p^{\ord_p(n)} } \sum_{f=1}^{ p^{\ord_p(n)} } N_{  \Z/ p^e\Z }^1(f) f^{w}   \\ 
&\overset{\eqref{r=1_density}}{=} \prod_{ p | n }  \frac{1}{ p^{\ord_p(n)} } 
\left[  
\left\{ \sum_{d=0}^{\ord_p(n)-1} (p^{\ord_p(n)- d } - p^{\ord_p(n)-d-1})   ( p^d)^{w} 
\right\}  
+ \left\{ 1\cdot  (p^{\ord_p(n)})^{w} 
\right\} \right]  \\
&=  \prod_{ p | n } 
\left[  
\left\{ ( 1 - p^{-1} ) \sum_{d=0}^{\ord_p(n)-1} \left( p^{w-1} 
\right)^d    \right\}  
+  p^{\ord_p(n)(w-1)} 
\right]
\\
&= 
\begin{cases}
\prod_{p|n} \left( ( 1 - p^{-1} ) \frac{ 1 - p^{(w
-1)\ord_p(n)} }{ 1 - p^{w
-1} } +  p^{\ord_p(n)(w 
-1)} 
\right)
\
&\text{if}\ w
\neq 1  \\
\prod_{p|n} \left(  1 + \left( 1 - \frac{1}{p} \right) \ord_p(n) \right)
\
&\text{if}\ w
= 1
\end{cases}
\\
&= 
\begin{cases}
\prod_{p|n} \left[ \left( \frac{ 1 - p^{-1} }{ 1 - p^{w-1} } \right)  +  p^{\ord_p(n)(w -1)}  
\left( 1 -  \frac{ 1 - p^{-1} }{ 1 - p^{w-1} }
\right) \right]
\
&\text{if}\ w
\neq 1  \\
\prod_{p|n} \left(  1 + \left( 1 - \frac{1}{p} \right) \ord_p(n) \right)
\
&\text{if}\ w
= 1
\end{cases}
\end{align*}
\qed

So, the case $r=1$ is a matter of triviality, and 
not much extra effort in considering general $w\in\C$.
However, the proof for general $r\in\N$ is more complicated.   
Observe for $0\leq d < e$, 
{\small
\begin{equation} \label{Easy_N}
\begin{split}
N_{ \Z/ p^e\Z }^r( p^d ) &=
\text{The coefficient of $x^d$ in}\\
&\quad \left\{ (p^e - p^{e-1})\cdot 1 +  (p^{e-1} - p^{e-2})\cdot x + \cdots + (p-1) x^{e-1} \right\}^r   \\
&= (p^e - p^{e-1})^r\times \text{The coefficient of $x^d$ in}\
\left\{ 1 + \left( \frac{x}{p} \right) + \cdots + \left( \frac{x}{p} \right)^{e-1} \right\}^r   \\
&= (p^e - p^{e-1})^r\times \text{The coefficient of $x^d$ in}\
\left\{ 1 - \left( \frac{x}{p} \right) \right\}^{-r}   \\
&= (p^e - p^{e-1})^r\times \text{The coefficient of $x^d$ in}\
\sum_{l=0}^{e-1} {}_rH_l \left( \frac{x}{p} \right)^l   \\
&= (p^e - p^{e-1})^r {}_rH_d\left( \frac{1}{p} \right)^d
\end{split}
\end{equation}
}

To evaluate $N_{ \Z/ p^e\Z }^r( p^e )$ and for further computation, 
we prove the following identities involving the repeated combinations,
which may be of independent interest:


\begin{lemma} \label{HI}
For $e,r \in \N$, set a degree $e-1$ polynomial of $x$ by
\begin{equation} \label{fdef}
f^r_e(x) := 
\sum_{k=0}^{e-1} {}_rH_k x^k
\end{equation}
Then,
\begin{subequations} 
\begin{gather}
f^r_e(1) = {}_{r+1}H_{e-1} = {}_{e+r-1}C_{e-1} = {}_{e+r-1}C_r = {}_eH_r    \label{at1}  \\
(1-x)^rf^r_e(y) + y^ef^e_r(1-x) =  (1-x)^r\Big( f^r_e(y) - y^ef^r_e(1) \Big) + y^ef^e_{r+1}(1-x) 
 \label{asymmetry}  \\
(1-x)^rf^r_e(x) + x^ef^e_r(1-x) = 1  \label{symmetry} 
\end{gather}
\end{subequations}
\end{lemma}

\begin{proof}
\eqref{at1} is an easy consequence of the Taylor expansions at $x=0$:
\begin{align*}
f^r_e(1)  &= \sum_{k=0}^{e-1} {}_rH_k = 
\text{The coefficient of $x^{e-1}$ in}\ \Big( \sum_{k=0}^{e-1} {}_rH_k x^k \Big)\Big( \sum_{k=0}^{e-1}  x^k \Big)  \\
&= \text{The coefficient of $x^{e-1}$ in}\ \Big(  ( 1 - x )^{-r}(1-x)^{-1} = (1-x)^{ -(r+1)}  \Big) \\
&=  {}_{r+1}H_{e-1} = {}_{e+r-1}C_{e-1} = {}_{e+r-1}C_r = {}_eH_r 
\end{align*}

Now \eqref{asymmetry} is an immediate consequence of \eqref{at1}:
\begin{align*}
&\quad (1-x)^rf^r_e(y) + y^ef^e_r(1-x) =  (1-x)^r\Big( f^r_e(y) - y^ef^r_e(1) \Big) + y^e\Big( f^e_{r}(1-x) + f^r_e(1) (1-x)^r \Big)  \\
&\overset{\eqref{fdef}\eqref{at1}}{=} (1-x)^r\Big( f^r_e(y) - y^ef^r_e(1) \Big) 
+ y^e\left( \left( \sum_{k=0}^{r-1} {}_eH_k(1-x)^k \right) + {}_eH_r(1-x)^r \right)
\\
&= (1-x)^r\Big( f^r_e(y) - y^ef^r_e(1) \Big) 
+ y^e \sum_{k=0}^{r} {}_eH_k(1-x)^k 
\overset{\eqref{fdef}}{=} (1-x)^r\Big( f^r_e(y) - y^ef^r_e(1) \Big) + y^ef^e_{r+1}(1-x) 
\end{align*}

To prove \eqref{symmetry}, consider the Taylor expansions of degree $e$ of $(1-x)^{-r}$ at $x=0$,
and degree $r$ of $x^{-e}= \left( 1 - (1-x) \right)^{-e}$ at $x=1$, respectively,
\begin{subequations} 
\begin{align}
f^r_e(x) &=  \sum_{k=0}^{e-1} {}_rH_k x^k = (1-x)^{-r} + x^e\cdot O( 1; x\to 0 )
\label{T0}
\\
f^e_r(1-x) &= \sum_{k=0}^{r-1} {}_eH_k (1-x)^k = x^{-e} + (1-x)^r\cdot O( 1; x\to 1 )  
\label{T1}
\end{align}
\end{subequations}
We now set the left hand side polynomial of $x$ in \eqref{fdef} to be
\begin{equation} \label{LHS}
l(x) := (1-x)^rf^r_e(x) + x^ef^e_r(1-x)
\end{equation}
Then, applying \eqref{T0} and \eqref{T1} respectively to \eqref{LHS}, we see
\begin{equation} \label{lat0}
\begin{split}
l(x) &= (1-x)^r\Big(  (1-x)^{-r} + x^e\cdot O( 1; x\to 0 ) \Big) + x^ef^e_r(1-x)  \\
&= 1 +  x^e\cdot \Big( (1-x)^rO( 1; x\to 0 ) +  f^e_r(1-x) \Big) 
\end{split}
\end{equation}
\begin{equation}  \label{lat1}
\begin{split} 
l(x) &= (1-x)^rf^r_e(x) + x^e \Big(  x^{-e} + (1-x)^r\cdot O( 1; x\to 1 )  \Big)   \\
&= 1 + (1-x)^r \cdot \Big(  f^r_e(x) + x^e\cdot O( 1; x\to 1 )  \Big)
\end{split}
\end{equation}
\eqref{lat0} and \eqref{lat1} respectively imply $l(x) - 1$ is divisible by
$x^e$ and $(1-x)^r$. Thus, $l(x) - 1$, a polynomial of $x$ of degree at most $e+r -1$,
is divisible by $x^e(1-x)^r$.  Of course, this implies $l(x) = 1$, and thus \eqref{symmetry}
has been proven.

\end{proof}




Now we may complete our computation of $N_{ \Z/ p^e\Z }^r( f )$:

\begin{proposition} \label{ Np }
For the cyclic ring
$\Z/ p^e\Z$ with $p$ a prime and $e\in \N$,
\begin{equation} \label{density}
N_{ \Z/ p^e\Z }^r( f )  =
\begin{cases}
(p^e - p^{e-1})^r {}_rH_d\left( \frac{1}{p} \right)^d\quad &\text{if}\ f = p^d\ (0\leq d < e)   \\
p^{e(r-1)} f^e_r\left( 1 - \frac{1}{p} \right)
= p^{e(r-1)}  \sum_{l=0}^{r-1}  {}_eH_l \left( 1 - \frac{1}{p} \right)^l  \quad &\text{if}\ f = p^e \\
0 \qquad &\text{if otherwise}
\end{cases}
\end{equation}
\end{proposition}

\begin{proof}[Proof of Proposition~\ref{ Np }]
It is obvious that $N_{ \Z/ p^e\Z }^r( f ) = 0$ if $f$ is not of the form $p^d\ (0\leq d\leq e)$.
The case of $0\leq d<e$ is already treated in \eqref{Easy_N}.


Finally, the case of $d=e$ can be taken care of as follows:
\begin{align*}
N_{ \Z/ p^e\Z }^r( p^e ) &=  \Big|  ( \Z/ p^e\Z )^r \Big| - \sum_{d=0}^{e-1} N_{ \Z/ p^e\Z }^r( p^d )   
\overset{\eqref{Easy_N}}{=} ( p^e )^r  - \sum_{d=0}^{e-1}   (p^e - p^{e-1})^r {}_rH_d\left( \frac{1}{p} \right)^d  \\
&\overset{\eqref{fdef}}{=}  ( p^e )^r  -   (p^e - p^{e-1})^r f^r_e\left( \frac{1}{p} \right)  
= p^{er}\Bigg( 1 - \left( 1 - \frac{1}{p} \right)^r  f^r_e\left( \frac{1}{p} \right)  \Bigg)
\\
&\overset{\eqref{symmetry}}{=} p^{er} \left( \frac{1}{p} \right)^ef^e_r\left( 1 - \frac{1}{p} \right)
= p^{e(r-1)} f^e_r\left( 1 - \frac{1}{p} \right)
= p^{e(r-1)}  \sum_{l=0}^{r-1}  {}_eH_l \left( 1 - \frac{1}{p} \right)^l \\
\end{align*}

\end{proof}

\begin{proof}[Proof of Main Theorem] 
{\small
\begin{align*}
&\quad
\frac{1}{n^r}\sum_{k_1, \ldots, k_r =1}^n \gcd(n,k_1\cdots k_r)^w
\overset{\eqref{reduction}}{=}
 \prod_{ p | n }  \frac{1}{ p^{r\ord_p(n)} } \sum_{f=1}^{ p^{\ord_p(n)} } N_{  \Z/ p^e\Z }^r(f) f^w \\
&\overset{\eqref{density}}{=} \prod_{ p | n }  \frac{1}{ p^{r\ord_p(n)} } 
\left[  
\left\{ \sum_{d=0}^{\ord_p(n)-1} (p^{\ord_p(n)} - p^{\ord_p(n)-1})^r {}_rH_d\left( \frac{1}{p} \right)^d  ( p^d)^w  \right\}  \right. \\
&\hspace{45mm} \left.
+ \left\{ p^{\ord_p(n)(r-1)}  \sum_{l=0}^{r-1}  {}_{\ord_p(n)}H_l \left( 1 - \frac{1}{p} \right)^l  \right\}  (p^{\ord_p(n)})^w \right]  \\
&=  \prod_{ p | n } 
\left[  
\left\{ ( 1 - p^{-1} )^r \sum_{d=0}^{\ord_p(n)-1}  {}_rH_d\left( p^{w-1} \right)^d    \right\}  
+ \left\{ p^{\ord_p(n)(w-1)}  \sum_{l=0}^{r-1}  {}_{\ord_p(n)}H_l \left( 1 - \frac{1}{p} \right)^l  \right\}   \right]
\\
&\overset{\eqref{fdef}}{=}  \prod_{ p | n } 
\Bigg[   ( 1 - p^{-1} )^r f^r_{\ord_p(n)}( p^{w-1} ) +  
\left( p^{w-1} \right)^{\ord_p(n)} f_r^{\ord_p(n)}   ( 1 - p^{-1} ) \Bigg]
\\
&\overset{\eqref{asymmetry}}{=}  
\begin{cases}
\prod_{ p | n } 
\Bigg[   ( 1 - p^{-1} )^r f^r_{\ord_p(n)}( p^{w-1} ) +  
\left( p^{w-1} \right)^{\ord_p(n)} f_r^{\ord_p(n)}   ( 1 - p^{-1} ) \Bigg]
\ 
&\text{if}\ w\neq 1
\\
\prod_{ p | n } 
\Bigg[   ( 1 - p^{-1} )^r \Bigg( f^r_{\ord_p(n)}( p^{w-1} ) -  \left( p^{w-1} \right)^{\ord_p(n)}f^r_{\ord_p(n)}(1) \Bigg)  +  
\left( p^{w-1} \right)^{\ord_p(n)} f_{r+1}^{\ord_p(n)}   ( 1 - p^{-1} ) \Bigg]
\
&\text{if}\ w = 1 
\end{cases}
\\
&\overset{\eqref{symmetry}}{=} 
\begin{cases}
\prod_{ p | n } 
\Bigg[   ( 1 - p^{-1} )^r 
\frac{ 1  -   ( p^{w-1
} )^{\ord_p(n)} f_r^{\ord_p(n)}( 1 - p^{w-1} } 
{  (   1  -   p^{w-1} )^r } 
+ 
\left( p^{w-1} 
\right)^{\ord_p(n)} f_{r}^{\ord_p(n)}   ( 1 - p^{-1} ) \Bigg]
\ 
&\text{if}\ w\neq 1 \\
\prod_{ p | n }  \Bigg[  f_{r+1}^{\ord_p(n)}   ( 1 - p^{-1} ) \Bigg]
\ 
&\text{if}\ w=1
\end{cases}
\\
&\overset{\eqref{fdef}}{=} 
\begin{cases}
\prod_{ p \mid n}
\left[
\left( \frac{ 1 - p^{-1} }{ 1 - p^{w-1}}
\right)^r
+ p^{\ord_p(n)(w 
-1)} \sum_{l = 0}^{r-1} {}_{ \ord_p(n) }H_l  
\left\{
\left( 1 - p^{-1} \right)^{l} 
- 
\left( \frac{ 1 - p^{-1} }{ 1 - p^{w-1}}
\right)^r
( 1 - p^{w-1}
)^{l}   \right\} 
\right]
\
&\text{if}\ w
\neq 1  \\
\prod_{ p \mid n}  
\left[
\sum_{l = 0}^{r} {}_{ \ord_p(n) }H_l  \left( 1 - \frac{1}{p} \right)^l 
\right]
\
&\text{if}\ w
= 1
\end{cases}
\\
&=
\begin{cases}
\prod_{ p \mid n}
\left[
\left( \frac{ 1 - p^{-1} }{ 1 - p^{w
-1} } \right)^r
+ p^{\ord_p(n)(w 
-1)} ( 1 - p^{-1} )^r \sum_{l = 0}^{r-1} {}_{ \ord_p(n) }H_l  
\left\{
\left( 1 - p^{-1} \right)^{l-r} 
- ( 1 - p^{w-1}
)^{l-r}   \right\} 
\right]
\
&\text{if}\ w
\neq 1  \\
\prod_{ p \mid n}  
\left[
\sum_{l = 0}^{r} {}_{ \ord_p(n) }H_l  \left( 1 - \frac{1}{p} \right)^l 
\right]
\
&\text{if}\ w
= 1
\end{cases}
\end{align*}
}
\end{proof}

\begin{proof}[Proof of Corollary\ref{asymptotic}]
When $w=1,$
\begin{multline*}
\lim_{r\to + \infty} \prod_{ p \mid n}  
\left[
\sum_{l = 0}^{r} {}_{ \ord_p(n) }H_l  \left( 1 - \frac{1}{p} \right)^l 
\right]
= \prod_{ p \mid n} \left[
\sum_{l = 0}^{+ \infty} {}_{ \ord_p(n) }H_l  \left( 1 - \frac{1}{p} \right)^l 
\right]
\\
\overset{ | 1 - \frac{1}{p} | < 1}{=} \prod_{ p \mid n}  \left( 1 -  \left( 1 - \frac{1}{p} \right) \right)^{- \ord_p(n)}
=  \prod_{ p \mid n}  p^{\ord_p(n)} = n
\end{multline*}

When $ | p^w - 1 | > 1$, since $| 1 - p^{-1} | < 1 <  | 1 - p^w |$,
{\small
\begin{align*}
&\ \lim_{r\to + \infty} 
\prod_{ p \mid n}
\left[
\left( \frac{ 1 - p^{-1} }{ 1 - p^{w-1}}
\right)^r
+ p^{\ord_p(n)(w 
-1)} \sum_{l = 0}^{r-1} {}_{ \ord_p(n) }H_l  
\left\{
\left( 1 - p^{-1} \right)^{l} 
- 
\left( \frac{ 1 - p^{-1} }{ 1 - p^{w-1}}
\right)^r
( 1 - p^{w-1}
)^{l}   \right\} 
\right]
\\
&= \prod_{ p \mid n}  \left[   \lim_{r\to + \infty} \left( \frac{ 1 - p^{-1} }{ 1 - p^{w-1}}
\right)^r
+  p^{\ord_p(n)(w 
-1)}   \sum_{l = 0}^{+\infty} {}_{ \ord_p(n) }H_l  
\left( 1 - p^{-1} \right)^{l}  \right.  \\
&\qquad\qquad \left.
- p^{\ord_p(n)(w 
-1)} \lim_{r\to + \infty}  ( 1 - p^{-1}
)^{r}  \sum_{l = 0}^{r-1} {}_{ \ord_p(n) }H_l ({ 1 - p^{w-1}}
)^{l-r}
\right]
\\
&\overset{
|  \frac{ 1 - p^{-1} }{ 1 - p^{w-1}} | <1,\ | 1 - p^{-1} | < 1
}{=} 
\prod_{ p \mid n}  \left[   
0
+  p^{\ord_p(n)(w 
-1)}   p^{\ord_p(n)}
\right.
\\
& \hspace{40mm}
\left.
- p^{\ord_p(n)(w 
-1)} \lim_{r\to + \infty}  ( 1 - p^{-1}
)^{r}  \sum_{l = 0}^{r-1} {}_{ \ord_p(n) }H_l  ( { 1 - p^{w-1}}
)^{l-r}
\right]
\\
&= \prod_{ p \mid n}  \left[   
p^{\ord_p(n)(w)} 
- p^{\ord_p(n)(w 
-1)} \lim_{r\to + \infty}  ( 1 - p^{-1}
)^{r}  \sum_{l = 0}^{r-1} {}_{ \ord_p(n) }H_l  ( { 1 - p^{w-1}}
)^{l-r}
\right]
\end{align*}
}
This is shown to be equal to $n^w$, if the following is shown:
\begin{equation*}
\lim_{r\to + \infty}  ( 1 - p^{-1}
)^{r}  \sum_{l = 0}^{r-1} {}_{ \ord_p(n) }H_l  ( { 1 - p^{w-1}}
)^{l-r}
= 0
\end{equation*}
However, this follows from
\begin{multline*}
\lim_{r\to + \infty}  \Big|  ( 1 - p^{-1}
)^{r}  \sum_{l = 0}^{r-1} {}_{ \ord_p(n) }H_l  ( { 1 - p^{w-1}}
)^{l-r}  \Big| 
=
\lim_{r\to + \infty} 
\frac{   \Big|   \sum_{l = 0}^{r-1} {}_{ \ord_p(n) }H_l  ( \frac{1}{ 1 - p^{w-1}}
)^{r-l}  \Big|  }
{ \big| (   1 - p^{-1} )^{-1} \big|^r }
\\
\leq 
\lim_{r\to + \infty} 
\frac{     \sum_{l = 0}^{r-1} {}_{ \ord_p(n) }H_l  \Big|  ( \frac{1}{ 1 - p^{w-1}}
)^{r-l}  \Big|  }
{ \big| (   1 - p^{-1} )^{-1} \big|^r }
\overset{  | p^w - 1 | > 1 }{\leq}
\lim_{r\to + \infty} 
\frac{     \sum_{l = 0}^{r-1} {}_{ \ord_p(n) }H_l  
\cdot 1
}
{ \big| (   1 - p^{-1} )^{-1} \big|^r }
\\
\overset{\eqref{fdef}}=  \lim_{r\to + \infty} \frac{  f^{\ord_p(n)}_r(1) }{  \big| (   1 - p^{-1} )^{-1} \big|^r } 
\overset{\eqref{at1}}=  \lim_{r\to + \infty} \frac{  {}_{\ord_p(n)+r-1}C_{\ord_p(n)} }{  \big| (   1 - p^{-1} )^{-1} \big|^r } 
= 0,
\end{multline*}
where the last equality follows since $ {}_{\ord_p(n)+r-1}C_{\ord_p(n)}$ is a degree $\ord_p(n)$ polynomial of $r$
and $\big| (   1 - p^{-1} )^{-1} \big| > 1$.
Now the claim has been proven.
\end{proof}


\end{document}